\documentclass[11pt,dvips]{amsart}  
\usepackage{amssymb,amsmath,amsfonts}  
\usepackage[dvips]{graphicx,xcolor}
\usepackage{cite}
\allowdisplaybreaks[4]
\definecolor{pgray}{gray}{0.8}

\newtheorem{theorem}{Theorem}[section]

\newtheorem{proposition}[theorem]{Proposition}

\numberwithin{equation}{section}

\title{\mbox{}}

\begin{document}
\begin{center}
{\bf \LARGE{
%
%
%
	Characteristics of drift effects arising from nonlinear symmetry of the quasi-geostrophic equation
}}\\
\vspace{5mm}
Masakazu Yamamoto\footnote{e-mail : \texttt{mk-yamamoto@gunma-u.ac.jp}}\\Graduate School of Science and Technology, Gunma University
\end{center}
\maketitle
\vspace{-15mm}
%
\begin{abstract}
This paper compares two similar diffusion equations that appear in meteorology.
One is the quasi-geostrophic equation, and the other is the convection-diffusion equation.
Both are two-dimensional bilinear equations, and the order of differentiation is the same.
Naturally, their scales also coincide.
However, the direction in which the nonlinear effects act differs: one acts along the isothermal surface, while the other acts along the temperature gradient in a specified direction.
The main assertion quantifies this difference through the large-time behavior of their solutions.
In particular, the nonlinear distortions in the asymptotic profiles of both equations are compared.
These profiles are uniquely determined based on the parabolic scale.
In this context, the spatial symmetry of the first approximation plays a crucial role, but the solutions require no symmetry.
As an appendix, the mixed problem of those models are studied.
\end{abstract}

\section{Introduction}
We treat mainly the Cauchy problem of the following quasi-geostrophic equation:
\begin{equation}\label{qg}
\left\{
\begin{array}{lr}
	\partial_t \theta + \mathcal{R}^\bot\theta\cdot\nabla \theta = \Delta \theta,
	&
	t > 0,~ x \in \mathbb{R}^2,\\
	\theta (0,x) = \theta_0 (x),
	&
	x \in \mathbb{R}^2,
\end{array}
\right.
\end{equation}
where $\partial_t = \partial/\partial t$ and $\mathcal{R}^\bot = (-\mathcal{R}^2, \mathcal{R}^1)$ for $\mathcal{R}^j \varphi = \partial_j (-\Delta)^{-1/2} \varphi = \mathcal{F}^{-1} [i\xi_j\mathcal{F} [\varphi]/|\xi|]$ is the `orthogonal' Riesz transform.
The unknown function $\theta = \theta (t,x) \in \mathbb{R}$ describes potential temperature.
This nonlinear term represents the effect of the Coriolis force $\mathcal{R}^\bot \theta$ acting along isothermal surfaces on the potential temperature.
For this problem, global in time solvability in Lebesgue spaces is guaranteed by Moser--Nash theory.
Moreover, the solution satisfies
\begin{equation}\label{decayth}
	\| \theta (t) \|_{L^q (\mathbb{R}^2)}
	\le
	C (1+t)^{-\gamma_q}
\end{equation}
and
\begin{equation}\label{decaythwt}
	\| x \theta (t) \|_{L^q (\mathbb{R}^2)}
	\le
	C t^{-\gamma_q} (1+t)^{1/2}
\end{equation}
for $1 \le q \le \infty$ and $\gamma_q = 1-\frac1q$ when $\theta_0 \in L^1 (\mathbb{R}^2) \cap L^\infty (\mathbb{R}^2)$ and $x\theta_0 \in L^1 (\mathbb{R}^2)$ are supposed.
The asymptotic profile of the solution as $t \to +\infty$ is given by
\begin{equation}\label{asymp0}
	\| \theta (t) - M_0 G(t) \|_{L^q (\mathbb{R}^2)} = o (t^{-\gamma_q}),
\end{equation}
where $M_0 = \int_{\mathbb{R}^2} \theta_0 (x) dx = \int_{\mathbb{R}^2} \theta (t,x) dx$ is the conserved quantity and $G(t) = (4\pi t)^{-1} e^{-|x|^2/4t}$ is the fundamental solution (see\cite{Schnbk2}).
Based on considerations of similar models, it is expected that the difference on \eqref{asymp0} contains nonlinear components.
We introduce the following convection-diffusion equation as a similar model:
\begin{equation}\label{bg}
\left\{
\begin{array}{lr}
	\partial_t \rho - \Delta \rho = a\cdot\nabla (|\rho|\rho),
	&
	t > 0,~ x \in \mathbb{R}^2,\\
	\rho (0,x) = \rho_0 (x),
	&
	x \in \mathbb{R}^2.
\end{array}
\right.
\end{equation}
The nonlinear force in this equation represents convection generated by temperature differences.
Both equations are bilinear, and the order of derivatives is the same if the Riesz transform is regarded as a zeroth-order derivative.
Here, note that $\mathcal{R}^\bot\theta\cdot\nabla \theta = \nabla \cdot (\theta \mathcal{R}^\bot \theta)$.
On the other hand, focusing on symmetry, the effects of these nonlinearities appear entirely different.
For $\rho$, global in time solvability and decay as $t \to +\infty$ as those of \eqref{qg} can be proven.
Escobedo and Zuazua \cite{Escbd-ZZ} discovered the logarithmic shift arising from nonlinearity.
The nonlinearity also provides the following distortion of solution symmetry.
\begin{proposition}\label{propBgs}
Let $\rho_0 \in L^1 (\mathbb{R}^2) \cap L^\infty (\mathbb{R}^2)$ and $x\rho_0 \in L^1 (\mathbb{R}^2)$.
Then the solution $\rho$ of \eqref{bg} fulfills that
\[
\begin{split}
	&\biggl\| \rho(t) - M_0 G(t) - \frac{|M_0| M_0}{8\pi} a\cdot \nabla G(t) \log t - M_1 \cdot \nabla G(t) - \frac{|M_0| M_0}{8\pi} \sum_{k=1}^\infty \left(\frac{t}2\right)^k \frac{a\cdot\nabla (-\Delta)^k G(t)}{k k!} \biggr\|_{L^q (\mathbb{R}^2)}\\
	&= o(t^{-\gamma_q-1/2})
\end{split}
\]
as $t \to +\infty$ for $1 \le q \le \infty$ and $\gamma_q = 1-\frac1q$, where $M_0 = \int_{\mathbb{R}^2} \rho_0 (x) dx \in \mathbb{R}$ and $M_1 = - \int_{\mathbb{R}^2} x\rho_0 (x)dx + a \int_0^\infty \int_{\mathbb{R}^2} ((|\rho|\rho)(s,y) - |M_0| M_0 G^2 (1+s,y)) dyds \in \mathbb{R}^2$.
In addition, if $|x|^2 \rho_0 \in L^1 (\mathbb{R}^2)$, then the left-hand side is estimated by $O (t^{-\gamma_q-1}(\log t)^2)$ as $t\to\infty$.
\end{proposition}
%
Here, the last term
\begin{equation}\label{distbg}
	J_1 (t) = \frac{|M_0| M_0}{8\pi} \sum_{k=1}^\infty \left(\frac{t}2\right)^k \frac{a\cdot\nabla (-\Delta)^k G(t)}{k k!}
\end{equation}
provides the nonlinear distortion.
Indeed, this function has the same scale as $M_1 \cdot \nabla G(t)$ that is
\begin{equation}\label{scJ}
	\lambda^{2+1} (M_1\cdot\nabla G,J_1)(\lambda^2 t, \lambda x) = (M_1\cdot\nabla G,J_1) (t,x)
\end{equation}
for $\lambda > 0$, and then
\[
	\| (M_1\cdot\nabla G,J_1) (t) \|_{L^q (\mathbb{R}^2)} = t^{-\gamma_q-\frac12} \| (M_1\cdot\nabla G,J_1) (1) \|_{L^q (\mathbb{R}^2)}
\]
for $t > 0$.
Both of $M_1 \cdot \nabla G$ and $J_1$ are odd-type functions.
However, any terms of $J_1$ exhibit higher symmetries than $M_1 \cdot \nabla G$.
Generally, in linear phenomena, based on the second law of thermodynamics, there is a clear correspondence between scale and symmetry.
Nonlinear distortions are what break this correspondence.
The convection diffusion is affected by the nonlinear distortion whenever $a \neq 0$ since $\{ a\cdot\nabla (-\Delta)^k G\}_{k\in\mathbb{N}}$ are independent.
Another point to note is that the coefficients in the asymptotic expansion contain the ambiguous component $\int_0^\infty \int_{\mathbb{R}^2} ((|\rho|\rho)(s,y) - |M_0| M_0 G^2 (1+s,y)) dyds$.
For the Burgers-type equation with the nonlinear term replaced by $a \cdot \nabla (\rho^2)$ instead of $a \cdot \nabla (|\rho| \rho)$, we see that
\[
	\rho (t) \sim M_0 G(t) + \frac{M_0^2}{8\pi} a \cdot \nabla G(t) \log t + M_1 \cdot \nabla G(t) + \frac{M_0^2}{8\pi} \sum_{k=1}^\infty \left(\frac{t}2\right)^k \frac{a\cdot\nabla (-\Delta)^k G(t)}{k k!}
\]
as $t \to +\infty$ for the same $M_0$ and $M_1 = - \int_{\mathbb{R}^2} x\rho_0 (x)dx + a \int_0^\infty \int_{\mathbb{R}^2} (\rho^2(s,y) - M_0^2 G^2 (1+s,y)) dyds$.

We examine similar considerations for the quasi-geostrophic equation.
This leads to two conflicting predictions as follows:
\begin{itemize}
\item
	The expansion for $\theta$ contains similar distortions as in Proposition \ref{propBgs} since \eqref{qg} resembles \eqref{bg}
\item
	The expansion for $\theta$ has a different form from one of Proposition \ref{propBgs} since the symmetry of nonlinear terms on \eqref{qg} and \eqref{bg} are entirely different
\end{itemize}
In fact, we establish our main results as follows.
\begin{theorem}\label{thm-qg}
Let $\theta_0 \in L^1 (\mathbb{R}^2) \cap L^\infty (\mathbb{R}^2)$ and $x \theta_0 \in L^1 (\mathbb{R}^2)$.
Then the solution $\theta$ of \eqref{qg} fulfills that
\[
	\| \theta (t) - M_0 G(t) - M_1 \cdot \nabla G(t) \|_{L^q (\mathbb{R}^2)}
	= o (t^{-\gamma_q-1/2})
\]
as $t \to +\infty$ for $1 \le q \le \infty$ and $\gamma_q = 1-\frac1q$, where $M_0 = \int_{\mathbb{R}^2} \theta_0 (x) dx \in \mathbb{R}$ and $M_1 = - \int_{\mathbb{R}^2} x \theta_0 (x) dx \in \mathbb{R}^2$.
In addition, if $|x|^2 \theta_0 \in L^1 (\mathbb{R}^2)$, then the left-hand side is estimated by $O (t^{-\gamma_q-1}\log t)$ as $t\to\infty$.
\end{theorem}
That is, the second prediction is true and no effect of the nonlinear drift is appearing in the profile.
Particularly, nonlinear distortions never be seen.
Furthermore, neither logarithmic shift nor coefficient ambiguity is observed.
The logarithmic growth in the sharp estimate is mitigated.
This result is natural from a meteorological perspective.
According to the linear diffution, the potential temperature asymptotically converges to the radially symmetric function as \eqref{asymp0}.
We note that the initial conditions $x \theta_0,~ |x|^2 \theta_0 \in L^1 (\mathbb{R}^2)$ eliminate troughs and ridges.
Therefore, the Coriolis force acts most strongly on the radially symmetric component.
This generates smooth rotation due to geostrophic winds, but such a rotation contributes little to the averaging of potential temperature.
Consequently, the nonlinear external forces do not affect the behavior of the solution.
What is important in \eqref{qg} is not that it contains the Riesz transform, but rather the direction in which it acts (see Appendix \ref{sectApp}).

The quasi-geostrophic equation was formulated as a dissipative model by Constantin, Majda and Tabak \cite{Cnstntn-Mjd-Tbk}.
Subsequently, one with $(-\Delta)^{\alpha/2}$ instead of $-\Delta$ was introduced.
In particular, the case $\alpha = 1$ is important as a critical phenomena in the abstract theory (cf.\cite{Cffrll-Vssr,Ch-Cnstntn-Wu,Cnstntn-Crdb-W,Ju,Kslv-Nzrv-Vlbrg,Mur}).
However, in the context of asymptotic analysis, there is no critical case in $\alpha \le 2$ since the scale critical is $\alpha = 3$ (see\cite{Cnstntn-W,crdb-crdb2,YmSg}).
As clearly expected, when $\alpha < 2$, the nonlinear effects become even weaker than those in our problem.
At first glance, our expansion contradicts to one in the preceding work.
In \cite{OhAdhIwbch}, the asymptotic expansion with higher-order that contains the component arising from the nonlinear term is derived.
Actually, this component decays faster than expected.
Although, this is still required to deal the higher-order expansion.

The convection-diffusion equation with $a \cdot \nabla (|\rho|^{p-1} \rho)$ for general $p$ was introduced as an indicator when considering the behavior of solutions to general equations of continuity.
This equation has the case $p = 3/2$ as critical.
For subcritical cases, the logarithmic shift in solutions is interested as a nonlinear effect (cf.\cite{BKL,DrCrpo,Escbd-ZZ,FkdSt,Ksb}).
We focus on symmetry distortion of the solution, which more clearly represents nonlinear characteristics in the subcritical case $p = 2$.

\vspace{3mm}

\noindent
\textbf{Notation.}
We often omit the spatial parameter from functions, for example, $u(t) = u(t,x)$.
In particular, $G(t) * u_0 = \int_{\mathbb{R}^2} G(t,x-y) u_0 (y) dy$ and $\int_0^t g(t-s) * f(s) ds = \int_0^t \int_{\mathbb{R}^2} g(t-s,x-y) f(s,y) dyds$.
We define the Fourier transform and its inverse by $\mathcal{F} [\varphi] (\xi) = (2\pi)^{-1}$ $\int_{\mathbb{R}^2} e^{-ix\cdot\xi} \varphi (x) dx$ and $\mathcal{F}^{-1} [\varphi] (x) = (2\pi)^{-1} \int_{\mathbb{R}^2} e^{ix\cdot\xi} \varphi (\xi) d\xi$, where $i = \sqrt{-1}$.
We denote the derivations by $\partial_t = \partial / \partial t,~ \partial_j = \partial / \partial x_j$ for $j = 1,2,~ \nabla = (\partial_1,\partial_2),~ \nabla^\bot = (-\partial_2,\partial_1)$ and $\Delta = \partial_1^2 + \partial_2^2$.
The Riesz transform is defined by $R_j \varphi = \partial_j (-\Delta)^{-1/2} \varphi = \mathcal{F}^{-1} [i\xi_j |\xi|^{-1} \mathcal{F} [\varphi]]$ for $j = 1,2$, and $\mathcal{R}^\bot = (-\mathcal{R}_2,\mathcal{R}_1)$ and $\mathcal{R} = (\mathcal{R}_1, \mathcal{R}_2)$.
For $\beta = (\beta_1,\beta_2) \in \mathbb{Z}_+^2 = (\mathbb{N} \cup \{ 0 \})^2,~ |\beta| = \beta_1 + \beta_2$.
The Lebesgue space and its norm are denoted by $L^q (\mathbb{R}^2)$ and $\| \cdot \|_{L^q (\mathbb{R}^2)}$, that is, $\| f \|_{L^q (\mathbb{R}^2)} = (\int_{\mathbb{R}^2} |f(x)|^q dx)^{1/q}$ for $1 \le q < \infty$ and $\| f \|_{L^\infty (\mathbb{R}^2)}$ is the essential supremum.
The heat kernel and its decay rate on $L^q (\mathbb{R}^2)$ are symbolized by $G(t,x) = (4\pi t)^{-1} e^{-|x|^2/(4t)}$ and $\gamma_q = 1-\frac1q$.
We employ Landau symbol.
Namely, $f(t) = o(t^{-\mu})$ and $g(t) = O(t^{-\mu})$ mean $t^\mu f(t) \to 0$ and $t^\mu g(t) \to c$ for some $c \in \mathbb{R}$ such as $t \to +\infty$ or $t \to +0$, respectively.
Various nonnegative constants are denoted by $C$.

\section{Proof of the main result}
We employ the Escobedo--Zuazua theory together with the renormalization.
This combination was developed to derive the large-time behavior of solutions to the Burgers equation and the Keller--Segel chemotaxis model (cf.\cite{KtM,Ymd}).
This is also useful to lead local structure of solutions to several diffusion equations (see\cite{IshgKwkm13}).
Here, we apply this tool to extract the nonlinear components of the solution.

\subsection{The EZ expansion with the renormalization}\label{subsectgen}
To treat both equations, we introduce the general integral equation as
\begin{equation}\label{ms}
	u(t) = e^{t\Delta} u_0 + \int_0^t \nabla e^{(t-s)\Delta} \cdot \boldsymbol{f} [u] (s) ds.
\end{equation}
Precisely, $(u_0,\boldsymbol{f}) = (\theta_0,-\theta \mathcal{R}^\bot \theta)$ for \eqref{qg} and $(u_0,\boldsymbol{f}) = (\rho_0,a|\rho|\rho)$ for \eqref{bg}.
Then \eqref{decayth} and \eqref{decaythwt} are guaranteed for $u$.
Hence, H\"older inequality and, in some case, the boundedness of $\mathcal{R}^\bot$ say that $\| \boldsymbol{f} [u] (t) \|_{L^p (\mathbb{R}^2)} \le C(1+t)^{-\gamma_p-1}$ for $1 \le p < \infty$ and $\gamma_p = 1 - \frac1p$.
Thus, applying $L^p$-$L^q$ estimate for the heat semigroup yields that
\[
	\biggl\| \int_0^t \nabla e^{(t-s)\Delta} \cdot \boldsymbol{f} [u] (s) ds \biggr\|_{L^q (\mathbb{R}^2)}
	= O (t^{-\gamma_q-1/2} \log t),
\]
and then
\begin{equation}\label{asymmp0}
	\| u(t) - M_0 G(t) \|_{L^q (\mathbb{R}^2)}
	= O (t^{-\gamma_q-1/2} \log t)
\end{equation}
as $t \to +\infty$ for $1 \le q \le \infty$ and $M_0 = \int_{\mathbb{R}^2} u_0 (x) dx$.
We expand $u$ based on the EZ theory with the renormalization.
Then,
\begin{equation}\label{bsp}
\begin{split}
	&u(t) = M_0 G(t) + \nabla G(t) \cdot \int_0^t \int_{\mathbb{R}^2} \boldsymbol{f}[M_0 G] (1+s,y) dyds  + M_1 \cdot \nabla G(t)
	+J_1 (t) + r_2 (t)
\end{split}
\end{equation}
for
\[
\begin{split}
	&M_1 = - \int_{\mathbb{R}^2} y u_0 (y) dy + \int_0^\infty \int_{\mathbb{R}^2} (\boldsymbol{f} [u] (s,y) - \boldsymbol{f}[M_0 G] (1+s,y)) dyds,\\
	&J_1 (t) = \int_0^t \int_{\mathbb{R}^2} (\nabla G(t-s,x-y) - \nabla G(t,x)) \cdot \boldsymbol{f} [M_0 G] (s,y) dyds
\end{split}
\]
and
\[
\begin{split}
	&r_2 (t) = \int_{\mathbb{R}^2} (G(t,x-y) - G(t,x) + y \cdot \nabla G(t,x)) u_0 (y) dy\\
	&+ \int_0^t \int_{\mathbb{R}^2} (\nabla G(t-s,x-y) - \nabla G(t,x)) \cdot (\boldsymbol{f} [u] - \boldsymbol{f} [M_0 G] )(s,y) dyds\\
	&- \nabla G(t) \cdot \int_t^\infty \int_{\mathbb{R}^2} (\boldsymbol{f} [u] (s,y) - \boldsymbol{f}[M_0 G] (1+s,y)) dyds.
\end{split}
\]
One may doupt that $\boldsymbol{f} [M_0 G] (s)$ in $J_1$ has a singularity as $s \to +0$.
In fact, the mean value theorem for $\nabla G$ and Lebesgue convergence theorem fill this gap.
The treatment for the first term of $r_2$ is well-known.
Combining \eqref{decayth}-\eqref{asymp0}
gives $\| x^\beta (\boldsymbol{f} [u] - \boldsymbol{f} [M_0 G]) \|_{L^p (\mathbb{R}^2)} = O(t^{-\gamma_p-\frac32+\frac{|\beta|}2}\log t)$ as $t \to +\infty$ for $1 \le p < \infty$ and $|\beta| \le 1$.
Therefore, the $L^p$-$L^q$ estimate with the mean value theorem concludes that
\[
	\| r_2 (t) \|_{L^q (\mathbb{R}^2)} = O (t^{-\gamma_q -1} (\log t)^2)
\]
as $t \to +\infty$.
\subsection{The case of convection-diffusion equation}
For \eqref{bg}, the second term of \eqref{bsp} provides the logarithmic shift.
Indeed, applying the scale says that
\[
	\int_0^t \int_{\mathbb{R}^2} \boldsymbol{f}[M_0 G] (1+s,y) dyds = \int_0^t (1+s)^{-1} ds \int_{\mathbb{R}^2} \boldsymbol{f}[M_0 G] (1,y) dy
\]
and
\[
	\int_{\mathbb{R}^2} \boldsymbol{f}[M_0 G] (1,y) dy = |M_0| M_0 a \int_{\mathbb{R}^2} G^2 (1,y) dy = \frac{ |M_0| M_0 a}{8\pi}.
\]
The elementary calculus yields that
\[
\begin{split}
	&J_1 (t) = |M_0| M_0 \int_0^t \int_{\mathbb{R}^2} (a \cdot \nabla G(t-s,x-y) - a \cdot \nabla G(t,x)) G^2 (s,y) dyds\\
	&= \frac{|M_0| M_0}{8\pi} \int_0^t s^{-1} (a \cdot \nabla G (t-\tfrac{s}2) - a \cdot \nabla G(t)) ds.
\end{split}
\]
Finally, Taylor theorem around $-\frac{s}2 = 0$ gives us \eqref{distbg}.
\subsection{The case of quasi-geostrophic equation}\label{subsectqg}
For \eqref{qg}, the logarithmic shift in \eqref{bsp} is vanishing since the integrand $\boldsymbol{f}[M_0 G] = - M_0^2 G\mathcal{R}^\bot G$ is odd-type.
For the same reason, $J_1$ is written as
\begin{equation}\label{Jqg}
	J_1 (t) = - M_0^2 \int_0^t \nabla e^{(t-s)\Delta} \cdot (G\mathcal{R}^\bot G) (s) ds.
\end{equation}
This is the effect of the Coriolis force on the radially symmetric potential $M_0 G$, and it is expected not to act on the original potential $\theta$.
Here, we put $g(t,r) = G (t,x)$ and $h(t,r) = (-\Delta)^{-1/2} G(t,x)$ for $r = |x|$, then $(G \mathcal{R}^\bot G)(t,x) = x^\bot \frac{(gh')(t,r)}{r}$, and we conclude that
\[
	J_1 (t) = 0
\]
since $\nabla \cdot (G\mathcal{R}^\bot G) = 0$ other than the origin.
The ambiguous component $\int_0^\infty \int_{\mathbb{R}^2} ((\theta\mathcal{R}^\bot \theta) (s,y) - M_0^2 (G\mathcal{R}^\bot G) (1+s,y)) dyds$ in $M_1$ is vanishing since $\mathcal{R}^\bot$ is skew-adjoint in $L^2 (\mathbb{R}^2)$.
The logarithmic growth on the sharp estimate is mitigated.
Indeed, since the logarithmic shift is omitted from \eqref{bsp}, we see that
\begin{equation}\label{asymp0cr}
	\| \theta (t) - M_0 G(t) \|_{L^q (\mathbb{R}^2)} = O (t^{-\gamma_q - \frac12})
\end{equation}
as $t \to + \infty$.
Employing this instead of \eqref{asymp0} yields that
\[
	\| r_2 (t) \|_{L^q (\mathbb{R}^2)} = O (t^{-\gamma_q-1} \log t)
\]
as $t \to +\infty$.

\appendix

\section{A force acting perpendicular to the isothermal surface and equal in effect to the Coriolis force}\label{sectApp}
Here, to see how uniquely the Coriolis force behaves in \eqref{qg}, we introduce the following artificial problem:
\begin{equation}\label{qgp}
\left\{
\begin{array}{lr}
	\partial_t \vartheta - \nabla \cdot (\vartheta\mathcal{R}\vartheta) = \Delta \vartheta,
	&
	t > 0,~ x \in \mathbb{R}^2,\\
	\vartheta (0,x) = \vartheta_0 (x),
	&
	x \in \mathbb{R}^2,
\end{array}
\right.
\end{equation}
where $\mathcal{R} = (\mathcal{R}^1, \mathcal{R}^2)$ is the `forward-directed' Riesz transform.
This equation exhibits a nonlinear external force equivalent to the Coriolis force.
However, unlike the quasi-geostrophic equation, this force acts perpendicular to the isothermal surface.
In this sense, \eqref{qgp} possesses a symmetry distinct also from that of the convection-diffusion equation.
Under the $L^1 (\mathbb{R}^2) \cap L^\infty (\mathbb{R}^2)$-framework, Moser--Nash theory guarantees global in time existence of solutions and \eqref{decayth} for \eqref{qgp}.
Indeed,
\[
	-\int_{\mathbb{R}^2} |\vartheta|^{q-2} \vartheta \nabla \cdot (\vartheta\mathcal{R}\vartheta) dx
	=
	\frac{q-1}q \int_{\mathbb{R}^2} \mathcal{R} \vartheta \cdot \nabla (|\vartheta|^q) dx
	=
	\frac{q-1}q \int_{\mathbb{R}^2} |\vartheta|^q \Lambda \vartheta dx
\]
for $\Lambda \vartheta = \mathcal{F}^{-1} [|\xi| \hat{\vartheta}]$ and, hence, Stroock--Varopoulos estimate treats this term (cf.\cite{crdb-crdb2}).
Now we recall Section \ref{subsectgen}.
The solution $\vartheta$ satisfies \eqref{ms} for $(u_0,\boldsymbol{f}) = (\vartheta_0,\vartheta\mathcal{R}\vartheta)$.
Then, the weighted estimate \eqref{decaythwt} is derived from the $L^p$-$L^q$ estimate together with \eqref{decayth}.
The logarithmic shift
\[
	\int_0^t \int_{\mathbb{R}^2} \boldsymbol{f}[M_0 G] (1+s,y) dyds = M_0^2 \int_0^t (1+s)^{-1} ds \int_{\mathbb{R}^2} (G\mathcal{R} G) (1,y) dy
\]
is vanishing since $G\mathcal{R} G$ is also odd-type.
Similarly, the unclear components within $M_1$ disappear.
On the other hand, the nonlinear distortion appears.
Indeed, on the same way as in Section \ref{subsectqg}, we see that
\[
\begin{split}
	&J_1 (t) = M_0^2 \int_0^t \nabla e^{(t-s)\Delta} \cdot (G\mathcal{R} G) (s) ds\\
	&= \frac{\sqrt{\pi}M_0^2}{2\pi} \int_0^t \int_0^\infty \int_0^\infty \lambda^{-3/2} e^{-\frac1{4\lambda}} \nabla G(t-s) * (G(s) \nabla G (s+\sigma^2\lambda)) d\lambda d\sigma ds\\
	&= \frac{\sqrt{\pi}M_0^2}{8\pi^2} \int_0^t \int_0^\infty \int_0^\infty \lambda^{-3/2} e^{-\frac1{4\lambda}} s (2s+\sigma^2\lambda)^{-2} \Delta G(t-\tfrac{s^2}{2s+\sigma^2 \lambda}) d\lambda d\sigma ds.
\end{split}
\]
Here, $t-\frac{s^2}{2s+\sigma^2 \lambda} \ge \frac{t}2$ for $(\lambda,\sigma,s) \in (0,\infty)^2 \times (0,t)$, then $\Delta G(t-\tfrac{s^2}{2s+\sigma^2 \lambda})$ is analytic around $-\tfrac{s^2}{2s+\sigma^2 \lambda} = 0$ and we see from the elemental calculus that
\[
\begin{split}
	&J_1 (t) = \frac{\sqrt{\pi}M_0^2}{8\pi^2} \sum_{k=0}^\infty \frac{(-\Delta)^{k+1} G (t)}{k!} \int_0^t \int_0^\infty \int_0^\infty
		\lambda^{-3/2} e^{-\frac1{4\lambda}} s^{-3} \left( \frac{s^2}{2s+\sigma^2 \lambda} \right)^{k+2}
	d\lambda d\sigma ds\\
	&= \frac{\sqrt{\pi}M_0^2}{4\pi} \sum_{k=0}^\infty \left( \frac{t}2 \right)^{k+\frac12} \frac{(2k-1)!!  (-\Delta)^{k+1} G (t)}{k!(2k+2)!!}.
\end{split}
\]
Remember that any terms in this function have the same scale, that is, \eqref{scJ} holds also for this $J_1$.
From the same reason as in Section \ref{subsectqg}, we see that $\| r_2 (t) \|_{L^q (\mathbb{R}^2)} = O (t^{-\gamma_q} \log t)$ as $t \to +\infty$.
Summarizing the above, we see the following result.
\begin{proposition}\label{prop-qgp}
Let $\vartheta_0 \in L^1 (\mathbb{R}^2) \cap L^\infty (\mathbb{R}^2)$ and $x \vartheta_0 \in L^1 (\mathbb{R}^2)$.
Then the solution $\vartheta$ of \eqref{qgp} fulfills that
\[
	\| \vartheta (t) - M_0 G(t) - M_1 \cdot \nabla G(t) - J_1 (t) \|_{L^q (\mathbb{R}^2)}
	= o (t^{-\gamma_q-1/2})
\]
as $t \to +\infty$ for the above $J_1$, and $1 \le q \le \infty$ and $\gamma_q = 1-\frac1q$, and $M_0 = \int_{\mathbb{R}^2} \vartheta_0 (x) dx \in \mathbb{R}$ and $M_1 = - \int_{\mathbb{R}^2} x \vartheta_0 (x) dx \in \mathbb{R}^2$.
In addition, if $|x|^2 \vartheta_0 \in L^1 (\mathbb{R}^2)$, then the left-hand side is estimated by $O (t^{-\gamma_q-1}\log t)$ as $t\to\infty$.
\end{proposition}
Compared to the previous results, unlike the convection-diffusion equation, no logarithmic shift nor ambiguous coefficients appears in the solution behavior.
From this perspective, the nonlinearity of \eqref{qgp} resembles that of the quasi-geostrophic equation.
On the other hand, unlike the quasi-geostrophic equation, the nonlinear distortion $J_1$ appears.
However, this distortion possesses a symmetry distinct from that of the convection-diffusion equation.
Precisely, $J_1$ for \eqref{qgp} is radially symmetric and the derivative order is shifted.
One final comment: like the convection-diffusion equation, this distortion reflects the nonlinear structure of \eqref{qgp}, i.e., $(-\Delta)^{k+1} G = - (-\Delta)^{k+1/2} \nabla \cdot \mathcal{R} G$.

\end{document}